\documentclass[11pt]{amsart}

\usepackage{times,epsf,amssymb,amsmath,rlepsf,hyperref}

\begin{document}

\newtheorem{thm}{Theorem}[section]
\newtheorem{lem}[thm]{Lemma}
\newtheorem{cor}[thm]{Corollary}

\theoremstyle{definition}
\newtheorem{defn}{Definition}[section]

\theoremstyle{remark}
\newtheorem{rmk}{Remark}[section]

\def\square{\hfill${\vcenter{\vbox{\hrule height.4pt \hbox{\vrule
width.4pt height7pt \kern7pt \vrule width.4pt} \hrule height.4pt}}}$}

\newenvironment{pf}{{\it Proof:}\quad}{\square \vskip 12pt}

\title{Foliations of Hyperbolic Space by Constant Mean Curvature Hypersurfaces}
\author{Baris Coskunuzer}
\address{Koc University \\ Department of Mathematics \\ Sariyer, Istanbul 34450 Turkey}
\email{bcoskunuzer@ku.edu.tr}
\thanks{The author is partially supported by EU-FP7 Grant IRG-226062, TUBITAK Grant 107T642 and TUBA-GEBIP Award.}

\maketitle


\newcommand{\cirD}{\overset{\circ}{D}}
\newcommand{\Si}{S^2_{\infty}({\mathbf H}^3)}
\newcommand{\SI}{S^n_{\infty}({\mathbf H}^{n+1})}
\newcommand{\PI}{\partial_{\infty}}
\newcommand{\BHH}{{\mathbf H}^{n+1}}
\newcommand{\BH}{{\mathbf H}^3}
\newcommand{\BR}{\mathbf R}
\newcommand{\BC}{\mathbf C}
\newcommand{\BZ}{\mathbf Z}

\begin{abstract}
We show that the constant mean curvature hypersurfaces in $\BHH$ spanning the boundary of a star
shaped $C^{1,1}$ domain in $\SI$ give a foliation of $\BHH$. We also show that if $\Gamma$ is a
closed codimension-$1$ $C^{2,\alpha}$ submanifold in $\SI$ bounding a unique constant mean
curvature hypersurface $\Sigma_H$ in $\BHH$ with $\PI \Sigma_H = \Gamma$ for any $H\in (-1,1)$,
then the constant mean curvature hypersurfaces $\{\Sigma_H\}$ foliates $\BHH$.
\end{abstract}

\section{Introduction}

In this paper, we are interested in the question of the existence of a foliation of $\BHH$ with constant mean curvature
(CMC) hypersurfaces asymptotic to the given codimension-$1$ submanifold in $\SI$. CMC hypersurfaces in hyperbolic space
became an area of active research after the progress on minimal hypersurfaces in hyperbolic space. In \cite{A1},
Anderson showed the existence of area minimizing hypersurfaces in $\BHH$ for any given codimension-$1$ closed
submanifold in $\SI$. He also showed that for any mean convex asymptotic boundary, there exists a unique minimal
hypersurface in hyperbolic space. Then, Hardt and Lin studied the regularity of these area minimizing hypersurfaces,
and generalized Anderson's uniqueness result to star-shaped domains in $\SI$ in \cite{HL}.

In the following decade, there have been many important generalizations of these results to CMC hypersurfaces in
hyperbolic space. In \cite{To}, Tonegawa generalized Anderson's existence, and Hardt and Lin's regularity results for
CMC hypersurfaces by using geometric measure theory methods. In the same year, by using similar techniques, Alencar and
Rosenberg got a similar existence result in \cite{AR}. By using analytical techniques, Nelli and Spruck generalized
Anderson's result on uniqueness of minimal hypersurfaces for mean convex domains to CMC context in \cite{NS}. Then,
Guan and Spruck extended Hardt and Lin's result of uniqueness of minimal hypersurfaces for star-shaped domains to CMC
hypersurfaces case in \cite{GS}.

After these uniqueness results, the question of how these CMC surfaces live together in $\BH$ became interesting.
Indeed, in \cite{CV}, Chopp and Velling studied this problem by using computational methods, and had an interesting
result that for many different type of curves in $\Si$, CMC surfaces give a foliation of $\BH$. However, they do not
classify any type of curves to give a foliation of $\BH$, but for some chosen special curves, they showed the existence
of a foliation with computer aid. In this paper, we show that for some specific classes of codimension-$1$ submanifolds
in $\SI$, the foliation of $\BHH$ by CMC surfaces exist.

\vspace{0.3cm}

\noindent \textbf{Theorem 3.3.} Let $\Gamma$ be the boundary of a star-shaped $C^{1,1}$ domain in $\SI$. Then, there
exists a foliation of $\BHH$ by CMC hypersurfaces $\{\Sigma_H\} $ where $\Sigma_H$ is a CMC hypersurface with mean
curvature $H \in (-1,1)$ and $\PI \Sigma_H = \Gamma$.

\vspace{0.3cm}

We also show that a similar result is true for any closed codimension-$1$ submanifold in $\SI$ bounding a unique CMC
hypersurface in $\BHH$ for any given $H \in (-1,1)$.

\vspace{0.3cm}

\noindent \textbf{Theorem 4.2.} Let $\Gamma$ be a $C^{2,\alpha}$ closed codimension-1 submanifold
in $\SI$. Also assume that for any $H\in (-1,1)$, there exists a unique CMC hypersurface $\Sigma_H$
with $\PI \Sigma_H =\Gamma$. Then, the collection of CMC hypersurfaces $\{\Sigma_H\}$ with $H\in
(-1,1)$ foliates $\BHH$.

\vspace{0.3cm}

The organization of the paper as follows. In Section 2, we give the basic definitions and results which we use
throughout the paper. In Section 3, we show the star-shaped case. In Section 4, we prove the general case. Finally in
Section 5, we have some concluding remarks.

\subsection{Acknowledgements:}

I am very grateful to the referee for very valuable comments and suggestions. I would like to thank Brian White, Yoshihiro Tonegawa
and Yair Minsky for very helpful comments.

\section{Preliminaries}

In this section, we will overview the basic results which we use in the following sections. Note that the results and
the notions in this section can also be found in the survey article \cite{Co3}. Let $\Sigma^n$ be a compact
hypersurface bounding a domain $\Omega^{n+1}$ in some ambient Riemannian manifold $X$. Let $A$ be the area of $\Sigma$,
and $V$ be the volume of $\Omega$. Let's vary $\Sigma$ through a one parameter family $\Sigma_t$, with the
corresponding area $A(t)$ and volume $V(t)$. If $f$ is the normal component of the variation, and $H$ is the mean
curvature of $\Sigma$, then we get $A'(0) = -\int_\Sigma n H f$, and $V'(0)=\int_\Sigma f$ where $n$ is the dimension
of $\Sigma$, and $H$ is the mean curvature.

Let $\Sigma$ be a hypersurface with boundary $\Gamma$ in $X$. Assume also that there exists a hypersurface $M$ with
$\partial M = \Gamma$ such that $M\cup\Sigma=\partial \Omega$ where $\Omega$ is a codimension-$0$ submanifold in $X$,
and mean curvature of $M$ is greater than $H_0>0$ everywhere (inward direction of $\Omega$). Define $V(t)$ to be the
volume of the domain $\Omega_t$ bounded by $M$ and $\Sigma_t$. Now, we define a new functional as a combination of $A$
and $V$. Let $I_H(t)= A(t) + n H V(t)$ for $H<H_0$. Note that $I_0(t)=A(t)$. If $\Sigma$ is a critical point of the
functional $I_H$ for any variation $f$, then this will imply that $\Sigma$ has constant mean curvature $H$ \cite{Gu},
\cite{AR}. Note that the critical points of the functional $I_H$ are independent of the choice of the hypersurface $M$
since if $\widehat{I}_H$ is the functional which is defined with a different hypersurface $\widehat{M}$, then $I_H -
\widehat{I}_H = C$ for some constant $C$. In particular, $H=0$ is the special case of minimal surfaces, for which the
theory is very well developed \cite{Ni}, \cite{CM}.

\begin{defn} Let $\Sigma$ be a hypersurface in a Riemannian manifold $X$. $\Sigma$ is called \textit{minimal hypersurface}
if it is critical point of $I_0$ (Area Functional) for any variation. Equivalently, $\Sigma$ has constant mean
curvature $0$ at every point. $\Sigma$ is called an \textit{area minimizing hypersurface} if $\Sigma$ is the absolute
minimum of the functional $I_0$ (having the smallest area) among hypersurfaces with the same boundary.
\end{defn}

For general $H$, the theory is similar, and called the constant mean curvature (CMC) case.

\begin{defn} Let $\Sigma$ be a hypersurface in a Riemannian manifold $X$. $\Sigma$ is called \textit{CMC hypersurface}
if it is a critical point of $I_H$ for any variation. Equivalently, $\Sigma$ has constant mean curvature $H$ at every
point. $\Sigma$ is a \textit{minimizing CMC hypersurface} if $\Sigma$ is the absolute minimum of the functional $I_H$
among hypersurfaces with the same boundary.
\end{defn}

\noindent \textbf{Notation: }From now on, we will call CMC hypersurfaces with mean curvature $H$ as
\textit{$H$-hypersurfaces}.

\vspace{0.3cm}

We will call any complete noncompact hypersurface $\Sigma_H$ as {\em minimizing $H$-hypersurface (area minimizing
hypersurface)} if any compact codimension-$0$ submanifold with boundary of $\Sigma_H$ is a minimizing $H$-hypersurface
(area minimizing hypersurface).

After these general definitions of minimal and $H$-hypersurfaces, we will quote some basic facts
about the $H$-hypersurfaces in hyperbolic space.

First, we fix a codimension-$1$ closed submanifold $\Gamma$ in $\SI$. $\Gamma$ separates $\SI$ into two parts, say
$\Omega^+$ and $\Omega^-$. By using these domains, we will give orientation to hypersurfaces in $\BHH$ asymptotic to
$\Gamma$. With this orientation, mean curvature $H$ is positive if the mean curvature vector points towards positive
side of the hypersurface, negative otherwise. The following fact is known as maximum principle.

\begin{lem} $[$Maximum Principle$]$
Let $\Sigma_1$ and $\Sigma_2$ be two hypersurfaces in a Riemannian manifold which intersect at a common point
tangentially. If $\Sigma_2$ lies in positive side (mean curvature vector direction) of $\Sigma_1$ around the common
point, then $H_1$ is less than or equal to $H_2$ ($H_1 \leq H_2$) where $H_i$ is the mean curvature of $\Sigma_i$ at
the common point. If they do not coincide in a neighborhood of the common point, then $H_1$ is strictly less than $H_2$
($H_1<H_2$).
\end{lem}

With a simple application of this maximum principle, by using horospheres in $\BHH$, it is easy to
show that if $\Sigma$ is a complete $H$-hypersurface in $\BHH$ asymptotic to a codimension-$1$
submanifold $\Gamma$ of $\SI$, then $|H|<1$ \cite{Co1}.

The following existence theorem for minimizing $H$-hypersurfaces in $\BHH$ asymptotic to $\Gamma$
for a given codimension-$1$ closed submanifold in $\SI$ was proved by Tonegawa \cite{To}, and
Alencar-Rosenberg \cite{AR} independently by using geometric measure theory methods.

\begin{thm} \cite{To}, \cite{AR}
Let $\Gamma$ be a codimension-$1$ closed submanifold in $\SI$, and let $|H|<1$. Then there exists a
minimizing $H$-hypersurface $\Sigma$ in $\BHH$ where $\PI \Sigma = \Gamma$. Moreover, any such
minimizing $H$-hypersurface is smooth except a closed singularity set of dimension at most $n-7$.
\end{thm}

Now, we define the convex hull of a subset $A$ of $\SI$ in $\BHH$. If $\gamma$ is a round $n-1$-sphere in $\SI$, then
there is a unique geodesic plane $P$ in $\BHH$ asymptotic to $\gamma$. $\gamma$ separates $\SI$ into two parts
$\Delta^+$ and $\Delta^-$. Similarly, $P$ divides $\BHH$ into two halfspaces $D^+$ and $D^-$ with $\PI D^\pm =
\Delta^\pm$. We will call the halfspace whose asymptotic boundary contains $A$ as \textit{supporting halfspace}. i.e.
if $A\subset\Delta^+$, then $D^+$ is a supporting halfspace.

\begin{defn} Let $A$ be a subset of $\SI$. Then the \textit{convex hull} of $A$,
$CH(A)$, is the smallest closed convex subset of $\BHH$ which is asymptotic to $A$. Equivalently,
$CH(A)$ can be defined as the intersection of all supporting closed half-spaces of $\BHH$
\cite{EM}.
\end{defn}

Note that the asymptotic boundary of the convex hull of a subset of $\SI$ is the subset itself,
i.e. $\PI(CH(A)) = A$. In general, we say the hypersurface $\Sigma$ has the convex hull property if
it is in the convex hull of its boundary, i.e. $\Sigma \subset CH(\partial\Sigma)$. In special
case, if $\Sigma$ is a complete and noncompact hypersurface in $\BHH$, then we say $\Sigma$ has
convex hull property if it is in the convex hull of its asymptotic boundary, i.e. $\Sigma\subset
CH(\PI\Sigma)$. The minimal hypersurfaces in $\BHH$ have convex hull property \cite{A1}. There is
also a generalization of this property to $H$-hypersurfaces in $\BHH$ \cite{Co1}.

Now, fix $\Gamma$ and orient all spheres accordingly. If $T$ is a round $n-1$-sphere in $\SI$, then there is a unique
$H$-hypersurface $P_H$ in $\BHH$ asymptotic to $T$ for $-1<H<1$ \cite{NS}. $T$ separates $\SI$ into two parts
$\Delta^+$ and $\Delta^-$. Similarly, $P_H$ divides $\BHH$ into two domains $D_H^+$ and $D_H^-$ with $\PI D_H^\pm =
\Delta^\pm$. We will call these regions as \textit{$H$-shifted halfspaces}. If the asymptotic boundary of a $H$-shifted
halfspace contains $\Gamma$, then we will call this $H$-shifted halfspace as \textit{supporting $H$-shifted halfspace}.
i.e. if $A\subset\Delta^+$, then $D_H^+$ is a supporting $H$-shifted halfspace.

\begin{defn} Let $\Gamma$ be a codimension-$1$ submanifold of $\SI$.
Then the {\em $H$-shifted convex hull} of $\Gamma$, $CH_H(\Gamma)$ is defined as the intersection
of all supporting closed $H$-shifted halfspaces of $\BHH$.
\end{defn}

Now, the generalization of convex hull property of minimal hypersurfaces in $\BHH$ to
$H$-hypersurfaces in $\BHH$ is as follows \cite{Co1}. Similar versions of this result have been
proved by Alencar-Rosenberg in \cite{AR}, and by Tonegawa in \cite{To}.

\begin{lem} \cite{To}, \cite{AR}, \cite{Co1}
Let $\Sigma$ be a $H$-hypersurface in $\BHH$ where $\PI\Sigma = \Gamma$ and $|H|<1$. Then $\Sigma$
is in the $H$-shifted convex hull of $\Gamma$, i.e. $\Sigma \subset CH_H(\Gamma)$.
\end{lem}

The next result which we give in this section is also a generalization of a result on area
minimizing hypersurfaces to $H$-hypersurfaces.

\begin{lem} \cite{Co1}
Let $\Gamma_1$ and $\Gamma_2$ be two disjoint codimension-$1$ closed submanifolds in $\SI$. If
$\Sigma_1$ and $\Sigma_2$ are minimizing $H$-hypersurfaces in $\BHH$ where $\PI \Sigma_i =
\Gamma_i$, then $\Sigma_1$ and $\Sigma_2$ are disjoint, too.
\end{lem}

Now, we give a regularity result for $H$-hypersurfaces near infinity due to Tonegawa.

\begin{lem} \cite{To}
Let $\Gamma$ be $C^{k,\alpha}$ codimension-1 submanifold in $\SI$ where $1\leq k  \leq n-1$ and
$0\leq \alpha \leq 1$ or $k=n$ and $0\leq \alpha <1$. If $\Sigma$ is a complete CMC hypersurface in
$\BHH$ with $\PI \Sigma = \Gamma$, then $\Sigma\cup\Gamma$ is a $C^{k,\alpha}$ submanifold with
boundary in $\overline{\BHH}$ near $\Gamma$.
\end{lem}

We finish this section with the strong maximum principle due to Simon \cite{Si}. The maximum principle we mentioned
above, Lemma 2.1, is true for smooth hypersurfaces. Since we are working with codimension-1 rectifiable currents
(Theorem 2.2), by the regularity results of the geometric measure theory, they might have $n-7$ dimensional singular
set. Hence, we need a stronger version of the maximum principle which applies to codimension-1 rectifiable currents.
Simon originally proved the following result for area minimizing rectifiable currents, and it naturally extends to
minimizing $H$-hypersurfaces, too. Other versions of strong maximum principle for more general settings can also be
found in \cite{SW} and \cite{Il}.

\begin{lem} \cite{Si} $[$Strong Maximum Principle$]$ Let $T_1$ and $T_2$ be two minimizing $H$-hypersurfaces in a
Riemannian manifold. Let $U$ be an open subset of $N$ such that $\partial T_1 =\partial T_2 = 0$ in $U$, and $\mbox{reg
} T_1 \ \cap \ \mbox{reg }T_2 \ \cap \  U =\emptyset$. Then $\mbox{spt }T_1\  \cap \ \mbox{spt }T_2 \  \cap \  U =
\emptyset$.
\end{lem}

In the proof of the main theorem of \cite{Si}, $T_i$ is described as the boundary of a codimension-$0$ rectifiable
current $E_i$. By defining $E_i$ as in the description of minimizing $H$-hypersurfaces in Section $1$ of \cite{AR} (or
Section $2$ of \cite{To}), the whole proof of \cite{Si} goes through for minimizing $H$-hypersurfaces case, too.

\section{CMC Foliation for Star Shaped Boundary}

In this section, we will show that if $\Gamma$ is the boundary of a star-shaped $C^{1,1}$ domain in $\SI$, then there
exists a foliation of $\BHH$ by CMC hypersurfaces $\{\Sigma_H\} $ where $\Sigma_H$ is a $H$-hypersurface with $H \in
(-1,1)$ and $\PI \Sigma_H = \Gamma$.

First, we need to show that the CMC hypersurfaces $\{\Sigma_H\}$ sharing ideal boundary $\Gamma$ are disjoint.

\begin{lem} Let $\Gamma$ be the boundary of a star-shaped domain in $\SI$. If $\Sigma_{H_1}$ and $\Sigma_{H_2}$ are minimizing
CMC hypersurfaces in $\BHH$ with $-1<H_1<H_2<1$, then $\Sigma_{H_1}$ and $\Sigma_{H_2}$ are disjoint.
\end{lem}

\begin{pf} Let's take upper half space model for $\BHH$. Let origin be the {\em center} point of the star shaped domain $\Omega^+$
which $\Gamma$ bounds in $\SI$. Let $\gamma$ be the unique geodesic in $\BHH$ connecting the origin and $\infty$-point
(the vertical line at $0$). Let ${\varphi_t(x)= tx}$ be the hyperbolic isometry which is a translation along the
geodesic $\gamma$. Then, $\Gamma_t=\varphi_t(\Gamma)$ be also a boundary of a star shaped domain in $\SI$. Note that as
$\Gamma$ is boundary of a star shaped domain, $\Gamma_t\cap\Gamma_s=\emptyset$ for any $t\neq s$. Hence, the family of
closed codimension-1 submanifolds $\{\Gamma_t\}$ in $\SI$ foliates $\SI- \{0,\infty\}$, i.e. $\SI - \{0,\infty\} =
\bigcup_{t\in(0,\infty)} \Gamma_t$.

Let $\Sigma_{H_i}$ be the minimizing $H_i$-hypersurface in $\BHH$ with $\PI \Sigma_{H_i} = \Gamma$. We suppose the mean
curvature vectors of $\Sigma_{H_1}$ and $\Sigma_{H_2}$ point into the domains of $H^{n+1}$ containing $\Omega_{+}$ (
$0\leq H_1<H_2$). The same argument adapts to the other cases. Let $\Sigma_{H_1}^t = \varphi_t(\Sigma_{H_1})$. Then
$\PI \Sigma_{H_1}^t = \Gamma_t$. By Lemma 2.4, $\Sigma_{H_1}^t\cap\Sigma_{H_1}^s=\emptyset$ where $t\neq s$. As
$\varphi_t$ is continuous family of isometries, this implies the family of minimizing $H_1$-hypersurfaces
$\{\Sigma_{H_1}^t\}$ foliates $\BHH$, i.e. $\BHH = \bigcup_{t\in(0,\infty)} \Sigma_{H_1}^t$.

If $\Sigma_{H_1}\cap \Sigma_{H_2} \neq \emptyset$, then let $t_0 = \inf_{t\in(1,\infty)} \{t \ : \Sigma_{H_1}^t \cap
\Sigma_{H_2} = \emptyset\}$. By $H$-shifted convex hull property (Lemma 2.3), there exists an infimum $1<t_0<\infty$.
This implies $\Sigma_{H_1}^{t_0}$ intersects $\Sigma_{H_2}$ tangentially and $\Sigma^{t_0}_{H_1}$ lies in positive side
of the $\Sigma_{H_2}$. However, by maximum principle (Lemma 2.1), this is a contradiction as $H_1<H_2$.
\end{pf}

The more general case of the lemma above will be proved in the following section. The following
result is essential for the proof of the main result of this section.

\begin{lem} \cite{GS} Let $\Gamma$ be the boundary of a star-shaped $C^{1,1}$ domain in $\SI$ and
$|H|<1$. Then there exists a unique $H$-hypersurface $\Sigma_H$ in $\BHH$ with $\PI \Sigma_H =
\Gamma$.
\end{lem}

Now, we show that the unique $H$-hypersurfaces in the above result indeed foliates $\BHH$.

\begin{thm}
Let $\Gamma$ be the boundary of a star-shaped $C^{1,1}$ domain in $\SI$ and $|H|<1$. Then there exists a foliation of
$\BHH$ by CMC hypersurfaces $\{\Sigma_H\} $ where $\Sigma_H$ is a $H$-hypersurface in $\BHH$ and $\PI \Sigma_H =
\Gamma$.
\end{thm}

\begin{pf} By Lemma 3.2, for any $H\in(-1,1)$, there  exists a minimizing $H$-hypersurface $\Sigma_H$ with $\PI \Sigma_H = \Gamma$.
By Lemma 3.1, for $H_1\neq H_2$, $\Sigma_{H_1}\cap\Sigma_{H_2} = \emptyset$. To prove the theorem, first we will show
that there is no gap between two $H$-hypersurface in the collection $\{\Sigma_H\}$. Then, we prove that
$\bigcup_{-1<H<1} \Sigma_H$ fills the whole $\BHH$. Observe that by Lemma 2.5, $\overline{\Sigma}_H$ is a $C^1$
submanifold in $\overline{\BHH}$, and hence $\Sigma_H$ is separating in $\BHH$. Let $\SI-\Gamma = \Omega^+ \cup
\Omega^-$ and $\BHH-\Sigma_H = D_H^+ \cup D_H^-$ with $\PI D_H^\pm = \Omega^\pm$.\\

\noindent {\bf Step 1:} There is no gap between two $H$-hypersurface in the collection $\{\Sigma_H\}$. i.e.
$D_{H_0}^-=\bigcup_{H<H_0} D_H^-$ and $D_{H_0}^+=\bigcup_{H>H_0} D_H^+$.\\

\begin{pf} Assume that there is a gap in the collection $\{\Sigma_H\}$. i.e. $D_{H_0}^- \supsetneq \bigcup_{H<H_0} D_H^-$ (The other case is similar).
Let $\{ H_i\ : \ H_i\in (H_0-\epsilon,H_0),\ i\in \BZ^+\}$ be an increasing sequence with $H_i\to H_0$. Our aim is to
build a sequence $\{S_i\}$ of compact $H_i$-hypersurfaces with $S_i \subset \Sigma_{H_i}$ such that a subsequence
$S_{i_j}$ converges to an $H_0$-hypersurface $\widehat{\Sigma}_{H_0}$ with $\PI \widehat{\Sigma}_{H_0} = \Gamma$. Then,
by using  $\widehat{\Sigma}_{H_0}$ and star-shapedness of $\Gamma$ (Lemma 3.2), we get a contradiction and prove that
there cannot be any gap between in the collection $\{\Sigma_H\}$.

Now, we construct the sequence $\{S_i\}$ of compact $H$-hypersurfaces with $S_i \subset \Sigma_{H_i}$. Let $x$ be a
point in the $H$-shifted convex hull of $\Gamma$, $CH_{H_0}(\Gamma)$. By Lemma 2.5, the regularity theorem for
$H$-hypersurfaces near infinity \cite{To}, for sufficiently small $\rho>0$, $\Sigma_H \cap \{x_{n+1} < \rho\}$ is a
graph over $\Gamma \times (0,\rho)$ in upper half space model for $\BHH$, i.e. $\Sigma_H \cap \{x_{n+1} < \rho\} =
u_H(\Gamma \times (0,\rho))$ for some function $u_H$. Let $R_0$ be sufficiently large so that $\partial (B_{R_0}(x)\cap
\Sigma_H) \subset \{x_{n+1} < \rho \}$ in upper half space model for $\BHH$ (Note that $\partial B_{R_0}(x)$ has mean
curvature $\coth{R_0}>1>H_0$). Let $\{R_i\}$ be an increasing sequence with $R_i>R_0$ for any $i$, and $R_i\to\infty$
as $i\to \infty$. Then, $\partial (B_{R_i}(x) \cap \Sigma_{H_i}) = \gamma_i = u_{H_i}(\alpha_i)$ where $\alpha_i$  is
an (n-1)-sphere in $\Gamma_0 \times (0,\rho)$. Hence, define $S_i = B_{R_i}(x) \cap \Sigma_{H_i}$. Then, $\partial S_i
= \gamma_i$ and $\gamma_i \to \Gamma$ as $i\to\infty$. Also, for some sufficiently large $K$, let
$N_K(CH_{H_0}(\Gamma))$ be the $K$-neighborhood of the $H_0$-shifted convex hull of $\Gamma$. Since for any $H_i$,
$\Sigma_{H_i} \subset CH_{H_i}(\Gamma)$, by replacing $\epsilon$ if necessary, we can assume that $S_i\subset
N_K(CH_{H_0}(\Gamma))$ for any $i$.

Now, for any $i$, $S_i$ is the minimizer of the functional $I_{H_i}(t)= A(t) + n H_i V(t)$ for fixed boundary $\gamma_i
\subset \partial B_{R_i}(x)$, where $A$ is the area of the codimension-1 submanifold in $B_{R_i}(x)$, and $V$ is the
volume of the component which the submanifold separates in $B_{R_i}(x)$ by \cite{To}, \cite{AR}. Observe that $H_i \to
H_0$ and $S_i \subset N_K(CH_{H_0}(\Gamma))$ with $\partial S_i=\gamma_i \to \Gamma$. By adapting the proof of theorem
3 in \cite{AR}, we get a subsequence $S_{i_j} \to \widehat{\Sigma}_{H_0}$ where $\widehat{\Sigma}_{H_0}$ is a
minimizing $H_0$-hypersurface with $\PI \widehat{\Sigma}_{H_0} = \Gamma$ as follows. Let $\Delta_n = B_n(x)$ where
$n\in \mathbf{N}$ with $n>N_0$ where $N_0$ is sufficiently large that $N_K(CH_{H_0}(\Gamma))$ separates $\Delta_{N_0}$
into two parts. We claim that for any $n>N_0$, we can find uniform area bounds $c_n,C_n>0$ such that $c_n< |S_i\cap
\Delta_n| <C_n$. Then by using the compactness theorem for integral currents \cite{Fe}, we get the desired limit
$\widehat{\Sigma}_{H_0}$. Fix $n>N_0$. For any $i>n$, $\gamma_i$ is not nullhomologous in $N_K(CH_{H_0}(\Gamma)) -
\Delta_n$ by construction. Hence, $S_i\cap \Delta_n \neq \emptyset$. However, since $S_i$ is a minimizing
$H_i$-hypersurface, $|S_i \cap \Delta_n| < |\partial \Delta_n|$ by Corollary 1.1 of \cite{AR}, and let $C_n=|\partial
\Delta_n|$.

Now, we define $c_n$. Let $\alpha$ be a round (n-1)-sphere in $\Omega^+$ and $P$ be the unique $H$-hypersurface which
is a hyperplane (an equidistant sphere) with $\PI P = \alpha$ and $P\cap N_K(CH_{H_0}(\Gamma)) =\emptyset$. Further
assume that, $P\cap \Delta_{N_0} \neq \emptyset$, by modifying $N_0$ if necessary. Recall that $N_K(CH_{H_0}(\Gamma))$
separates $\Delta_n$ into two components. Fix a point $q_n$ in the component not intersecting $P$. For any point $p$ in
$P\cap \Delta_n$, let $l_p$ be the unique geodesic segment connecting $p$ and $q_n$. We claim that any geodesic segment
$l_p$ starting at $p\in P\cap \Delta_n$ intersects $S_i\cap \Delta_n$ for any $i$. We can see this as follows.
$N_K(CH_{H_0}(\Gamma))$ separates $\Delta_n$ into two parts. Any $l_p$ connects these two parts in $\Delta_n$ (Note
that since $\Delta_n$ is convex, $l_p\subset \Delta_n$). Since $l_p$ connects different components of $\BHH-
N_K(CH_{H_0}(\Gamma))$ in $\Delta_n$, and $\gamma_i$ is not nullhomologous in $N_K(CH_{H_0}(\Gamma)) - \Delta_n$,
$\gamma_i$ cannot be nullhomologous in $N_K(CH_{H_0}(\Gamma)) - l_p$ either. Hence, for any $i$, $S_i\cap l_p\neq
\emptyset$. Since $\Delta_n \cap N_K(CH_{H_0}(\Gamma))$ is compact, the geodesic projection of $S_i\cap \Delta_n$ from
$q_n$ onto $P$ has bounded distortion. Then, there exists $\lambda_n>0$ such that $|S_i\cap \Delta_n| \geq \lambda_n
|P\cap \Delta_n|$ and define $c_n = \lambda_n |P\cap \Delta_n|$. Hence, we found uniform area bounds $c_n,C_n>0$ such
that $c_n< |S_i\cap \Delta_n| <C_n$. Then by using the compactness theorem for integral currents (\cite{Fe}, Theorem
4.2.17) and isoperimetric inequality, for each fixed $n$, a subsequence $S_{i_j}\cap \Delta_n$ converges to a
minimizing $H_0$-minimizing hypersurface $\widehat{\Sigma}^n_{H_0}$ in $\Delta_n$ ( by the lower bound $c_n$,
$\widehat{\Sigma}^n_{H_0}$ is nonempty for each $n$). By repeating the argument for each $n>N_0$ for the new sequence,
the diagonal sequence argument gives a subsequence $S_{i_j}$ converging to an integral $n$-current
$\widehat{\Sigma}_{H_0}$ on any compact set, in the weak topology. By construction, $\widehat{\Sigma}_{H_0}$ is a
minimizing $H_0$-hypersurface in $\BHH$ with $\PI \widehat{\Sigma}_{H_0} = \Gamma$. Note also that by construction,
$\widehat{\Sigma}_{H_0}$ may not be smooth everywhere, and it might have at most $n-7$ dimensional singular set by the
regularity theorems of the geometric measure theory \cite{Fe}.

Now, because of the gap, $\widehat{\Sigma}_{H_0}$ cannot be same with $\Sigma_{H_0}$. Hence, we get two different
minimizing $H_0$-hypersurface with asymptotic boundary $\Gamma$. We cannot directly use the uniqueness statement in
Lemma 3.2, and get the contradiction because the uniqueness statement is for the smooth hypersurfaces while
$\widehat{\Sigma}_{H_0}$ may not be smooth everywhere. Hence, we get the contradiction as follows. We know that
$\Gamma$ is star shaped. Let $\varphi_t$, $\Gamma_t$, and $\Sigma_{H_0}^t=\varphi_t(\Sigma_{H_0})$ be as in the proof
of previous lemma. Since $\Gamma_t\cap\Gamma_s=\emptyset$ when $t\neq s$, by Lemma 2.4 and Lemma 2.6,
$\widehat{\Sigma}_{H_0}\cap \Sigma_{H_0}^t = \emptyset$ for any $t\neq 1$. Since $\Sigma_{H_0}^t$ foliates whole $\BHH$
by construction, this implies $\widehat{\Sigma}_{H_0} = \Sigma_{H_0}^1=\Sigma_{H_0}$, which is a contradiction. The
other cases are similar. Hence, this shows that there is no gap between two $H$-hypersurface in the collection
$\{\Sigma_H\}$.
\end{pf}

\noindent {\bf Step 2:} The collection $\{\Sigma_H\}$ fills $\BHH$, i.e. $\BHH=\bigcup_{-1<H<1} \Sigma_H$.\\

\begin{pf} Assume that $\BHH\neq\bigcup_{-1<H<1} \Sigma_H$. Recall that $\SI-\Gamma = \Omega^+ \cup \Omega^-$ and
$\BHH-\Sigma_H = D_H^+ \cup D_H^-$ with $\PI D_H^\pm = \Omega^\pm$. If $\BHH\neq\bigcup_{-1<H<1} \Sigma_H$, then either
$\bigcap_{-1<H<1} D^+_H \neq \emptyset$ or $\bigcap_{-1<H<1} D^-_H \neq \emptyset$. Note that $D^+_{H_1}\varsubsetneq
D^+_{H_2}$ when $H_1<H_2$, and likewise $D^-_{H_1}\varsubsetneq D^-_{H_2}$ when $H_1>H_2$. Without loss of generality,
assume that $Z=\bigcap_{-1<H<1} D^-_H \neq \emptyset$. Now, the idea is to construct a minimizing $1$-hypersurface
$\widehat{\Sigma}_{1}$ with $\PI \widehat{\Sigma}_{1} = \Gamma$, by using the construction in previous part. In the
previous construction, we used $N_K(CH_{H_0}(\Gamma))$ as barrier, however we cannot define $1$-shifted convex hull
$CH_1(\Gamma)$ by definition. Instead of $CH_1(\Gamma)$, we define $\widehat{C}= N_K(X)$ where $X= D^-_{1-\epsilon} -
Z$. Then for any $H_i\subset (1-\epsilon, 1)$, $\Sigma_{H_i}\subset \widehat{C}$. Let $H_i\subset (1-\epsilon, 1)$ and
$H_i\nearrow 1$, and define the sequence $\{S_i\}$ of compact $H_i$-hypersurfaces with $S_i = B_{R_i}(x) \cap
\Sigma_{H_i}$ as before. Hence, we can get the upper bounds $C_n$ such that  $|\Delta_n \cap S_i|< C_n=|\partial
B_n(x)|$ as before. For lower bound $c_n$, since we cannot define $1$-hypersurface $P$ with $\PI P = \alpha$ in that
construction, we use ($1-\epsilon$)-hypersurface $P'$ with $\PI P' = \alpha$ to define the lower area bounds $c_n$ to
ensure that the limit is nonempty. An alternative way to show that we get nonempty limit is as follows. Let $z$ be a
point in $Z$, let $y$ be a point in $\Sigma_{1-\epsilon}$, and let $\beta$ be a path between $z$ and $y$. Since
$D^-_{H_1}\varsubsetneq D^-_{H_2}$ when $H_1>H_2$, for any $i$, $\beta\cap S_i\neq \emptyset$ and let $w_i$ be a point
in $\beta \cap S_i$. As $\beta$ compact, any subsequence of $\{w_i\}$ will have a limit point in $\BHH$. Hence
$\{S_i\}$ cannot escape to the infinity, and we have a nonempty limit. Hence, the proof in previous part goes through,
and a subsequence $S_{i_j}$ converges to a minimizing $1$-hypersurface $\widehat{\Sigma}_{1}$ with $\PI
\widehat{\Sigma}_{1} = \Gamma$.

Now, this gives a contradiction as there cannot be a $1$-hypersurface in $\BHH$ because of the horospheres ($\pm
1$-hypersurfaces), and the maximum principle. In particular, let $p$ be a point in  $\Omega^+\subset \SI$. Let $H_R^p$
be the horospheres in $\BHH$ with $\PI H_R^p = \{p\}$ and $R$ is the radius of the horosphere in Euclidean metric on
$\overline{\BHH}$ in Poincare ball model (Hence, $0<R<1$). Let $c=\sup \{ 0<R<1 \ | \
H^p_R\cap\widehat{\Sigma}_{1}=\emptyset \ \}$. Then, $H^p_c$ is a minimizing $1$-hypersurface which intersects
$\widehat{\Sigma}_{1}$ at a point $q$ (the point of first touch) and completely lie in one side of
$\widehat{\Sigma}_{1}$. If $q$ is a  regular (smooth) point of $\widehat{\Sigma}_{1}$, this is a contradiction by the
maximum principle (Lemma 2.1). If $q$ is not a regular point of $\widehat{\Sigma}_{1}$ then the contradiction comes
from the strict maximum principle in \cite{Si} (Lemma 2.6).
\end{pf}

By Step 1 and Step 2, the collection of minimizing $H$-hypersurfaces $\{\Sigma_H\}$ gives a foliation of $\BHH$, hence
the proof follows.
\end{pf}

\section{CMC Foliation for General Case}

In this section, we will show that if $\Gamma$ is a $C^{2,\alpha}$ closed codimension-1 submanifold
in $\SI$ such that for any $H\in (-1,1)$, there exists a unique CMC hypersurface $\Sigma_H$ with
$\PI \Sigma_H =\Gamma$, then the collection of CMC hypersurfaces $\{\Sigma_H\}$ foliates $\BHH$. In
other words, in this section, we will show that the star-shapedness condition on $\Gamma$ in
previous section is not essential, but bounding a unique CMC hypersurface for any $H\in (-1,1)$ is
the key condition.

First, we will generalize Lemma 3.1, the pairwise disjointness result in previous section to a new setting.

\begin{lem} Let $\Gamma$ be a $C^{2,\alpha}$ closed codimension-1 submanifold
in $\SI$. If $\Sigma_{H_1}$ and $\Sigma_{H_2}$ are minimizing CMC hypersurfaces in $\BHH$ with $\PI
\Sigma_{H_i}=\Gamma_i$ and $-1<H_1<H_2<1$, then $\Sigma_{H_1}$ and $\Sigma_{H_2}$ are disjoint.
\end{lem}

\begin{pf}
By Lemma 2.5, the regularity theorem for $H$-hypersurfaces near infinity \cite{To}, for sufficiently small $\rho>0$,
$\Sigma_H \cap \{x_{n+1} < \rho\}$ is a graph over $\Gamma \times (0,\rho)$ in upper half space model for $\BHH$, i.e.
$\Sigma_H \cap \{x_{n+1} < \rho\} = u_H(\Gamma \times (0,\rho))$. Let $p\in \BHH$ be a point, and $R_0>0$ be
sufficiently large that $\partial B_{R_0}(p)=S_{R_0}(p) \subset \{x_{n+1}<\rho\}$ (Note that $S_{R_0}(p)$ has mean
curvature $\coth{R_0}>1>H_2$). By p.599 in \cite{To}, the estimate $u_y=\frac{H}{\sqrt{1-H^2}}
\sqrt{1+|\nabla\varphi(x)|^2}$ implies that for sufficiently small $\epsilon>0$, $u_{H_1} (x,\epsilon) < u_{H_2}
(x,\epsilon)$ where $x\in \Gamma$ and $0\leq H_1<H_2$. Hence, if $\gamma_i = S_{R_0}(p) \cap \Sigma_{H_i}$, then
$\gamma_1\cap\gamma_2 = \emptyset$.

Let $S_i = \Sigma_{H_1}\cap B_{R_0}(p)$ and $\gamma_i = \partial S_i$. Also, let $\Omega_i$ be the region which $S_i$
separates in $B_{R_0}(p)$ where the mean curvature vector  on $S_i$ points outside of $\Omega_i$ (without loss of
generality, we assume $0<H_1<H_2<1$. If both negative, just change the direction. If different signs, same argument
works as area minimizing hypersurface $\Sigma_0$ would be a barrier between them \cite{AR}). An alternative way to
think about this picture is the Poincare ball model in $\BHH$ (See Figure 1). Hence, if we can show that $\Omega_2$ is
strictly included in $\Omega_1$, we are done ($S_1$ cannot touch $S_2$ because of maximum principle, Lemma 2.1). This
is because this implies $S_1 \cap S_2 = \emptyset$, and we already know $\Sigma_{H_1} - S_1$ is disjoint from
$\Sigma_{H_2} - S_2$ by \cite{To}.

\begin{figure}[t]

\relabelbox  {\epsfxsize=3in

\centerline{\epsfbox{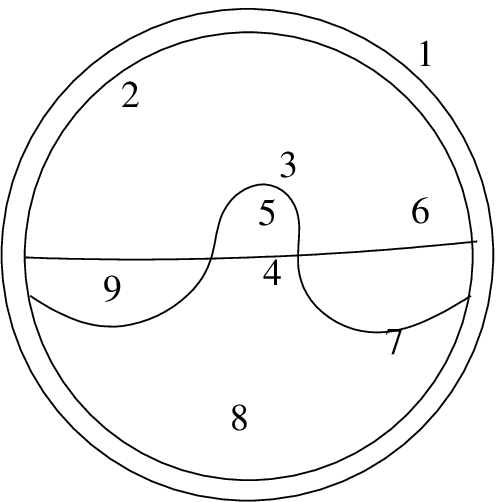}}}

\relabel{1}{$\SI$}

\relabel{2}{$S_{R_0}$}

\relabel{3}{$T_2$}

\relabel{4}{$T_1$}

\relabel{5}{$Q$}

\relabel{6}{$S_1$}

\relabel{7}{$S_2$}

\relabel{8}{$\Omega_1$}

\relabel{9}{$\Omega_2$}

\endrelabelbox

\caption{\label{fig:figure1} \small For $0<H_1<H_2<1$, $S_1$ is above $S_2$ near the boundary of
the ball $B_{R_0}(p)$ by \cite{To}}.

\end{figure}

Assume that $\Omega_2-\Omega_1 \neq \emptyset$. Since $\Sigma_{H_i}$ are minimizing
${H_i}$-hypersurfaces, then so is $S_i$, and hence $|S_i| + n H_i |\Omega_i|$ are minimizers for
the functional $I_{H_i}(t)= A(t) + n H_i V(t)$ for the fixed boundary $\gamma_i$. Consider
$\Omega_2'=\Omega_2 \cap \Omega_1$  and $S_2 ' = \partial \Omega_2' - S_{R_0}(p)$ which is another
hypersurface in $B_{R_0}(p)$ with $\partial S_2'=\gamma_2$. Let $Q = \Omega_2 - \Omega_1$. Let $T_1
=\partial Q \cap S_1$, $T_2 = \partial Q \cap S_2$ with $\widehat{S}_2 = S_2\cap \Omega_1$. Then,
$S_2 = \widehat{S}_2 \cup T_2$, $S_2' = \widehat{S}_2\cup T_1$, and $\Omega_2 = \Omega_2' \cup Q$.
Since $S_2$ is a minimizing $H_2$-hypersurface, $|S_2| + n H_2 |\Omega_2| \leq |S_2'| + n H_2
|\Omega_2'|$. This implies $(|\widehat{S}_2|+ |T_2|) + n H_2 (|\Omega_2'|+|Q|) \leq
(|\widehat{S}_2|+|T_1|) + n H_2 |\Omega_2'|$. After cancelations, we get $|T_2| + n H_2 |Q| \leq
|T_1|$ .

Similarly, let $\widehat{S}_1=S_1-T_1$ and $S_1' = \widehat{S}_1 \cup T_2$ and $\Omega_1' =
\Omega_1 \cup Q$. Again, since $S_1$ is minimizing $H_1$-hypersurface, then $|S_1| + n H_1
|\Omega_1| \leq |S_1'| + n H_1 |\Omega_1'|$. This implies $(|\widehat{S}_1|+ |T_1|) + n H_1
(|\Omega_1|) \leq (|\widehat{S}_1|+|T_2|) + n H_1 (|\Omega_1|+|Q|)$. After cancelations, we get
$|T_1| \leq |T_2| + n H_1 |Q|$. If we combine this with the previous inequality, we get $|T_2| + n
H_2 |Q| \leq |T_1|\leq |T_2| + n H_1 |Q|$. This implies $H_2\leq H_1$ which is a contradiction.
\end{pf}

In a sense, this lemma is a generalization of Lemma 3.1. Now, we generalize Theorem 3.3. The technique is basically same.

\begin{thm}
Let $\Gamma$ be a $C^{2,\alpha}$ closed codimension-1 submanifold in $\SI$. If for any
$H\in(-1,1)$, there exists a unique CMC hypersurface $\Sigma_H$ with $\PI \Sigma_H =\Gamma$, then,
the collection of CMC hypersurfaces $\{\Sigma_H\}$ with $H\in (-1,1)$ foliates $\BHH$.
\end{thm}

\begin{pf} By Lemma 4.1, $\Sigma_{H_1}\cap\Sigma_{H_2} = \emptyset$ for $H_1\neq H_2$. By using Lemma 4.1 instead of Lemma 3.1,
the proof of Step 1 of Theorem 3.3 adapts to this case easily, hence there cannot be any gap in the collection
$\{\Sigma_H\}$. Similarly, the proof of Step 2 of Theorem 3 also goes through in this case, and this shows
$\{\Sigma_H\}$ fill out the whole $\BHH$. The proof follows.
\end{pf}

\section{Final Remarks}

First we should note that Nelli-Spruck's uniqueness result does not automatically give a foliation
when it is combined with our result.

\begin{thm} \cite{NS} Let $\Omega$ be a $C^{2,\alpha}$ mean convex domain in $\SI$ and $\Gamma=\partial \Omega$.
Then for each $H\in(0,1)$ there exists a complete embedded hypersurface $M$ of $\BHH$ of constant mean curvature
$H$ with $\PI M = \Gamma$. Moreover, $M$ can be represented as a graph $x_{n+1}=u(x)$ over $\Omega$ with
$u\in C^{2,\alpha}(\overline{\Omega})$ and there is a unique such a graph.
\end{thm}

Unfortunately, this result for $C^{2,\alpha}$ mean convex domains does not give uniqueness in general. Other than the
obvious reason $H\in(0,1)$, it only asserts the uniqueness of the graphs, where there still might be other CMC
hypersurfaces which may not be represented as graph. So, for this case, Theorem 4.2 cannot be applied directly.
However, for the question whether the collection of CMC hypersurfaces $\{\Sigma_H\}_{H_\in(0,1)}$ in the theorem above
is a foliation of graphs or not, our techniques might be used, if one can show that the pairwise disjointness of the
graphs, and that the limit of the graphs is also a graph in that setting. Hence, with pairwise disjointness, the only
possibility not to have a foliation would be to have a gap, but the existence of a gap in the foliation would
contradict the uniqueness of graphs in the theorem above.

For the nonexistence question, by using the techniques of \cite{Co2}, it is possible to show the existence of a simple
closed curve $\Gamma$ in $\Si$ such that there is no foliation of $\BH$ by CMC hypersurfaces $\{\Sigma_H\}$ with $\PI
\Sigma_H = \Gamma$. In particular, such a foliation implies the uniqueness of not only minimizing CMC hypersurfaces,
but also any type of CMC hypersurfaces with asymptotic boundary $\Gamma$. In dimension $3$, by Hass's result in
\cite{Ha} and Anderson's result in \cite{A2}, there are examples of Jordan curves in $\Si$ bounding many minimal
surfaces in $\BH$. Hence, for those curves in $\Si$, there cannot be a foliation of $\BH$ by CMC surfaces because of
the maximum principle.

On the other hand, recently in \cite{Wa}, Wang showed that if a quasi-Fuchsian $3$-manifold $M$ contains a minimal
surface whose principle curvature is less than $1$, than $M$ admits a foliation by CMC surfaces by using volume
preserving mean curvature flow. If we lift this foliation to the universal cover, we get a foliation of $\BH$ by CMC
surfaces with same asymptotic boundary which is the limit set of $M$. However, the limit set of quasi-Fuchsian
manifolds are far from being smooth, even they contain no rectifiable arcs (\cite{Be}). Existence of one smooth point
in the limit set implies the group being Fuchsian, which means the limit set is a round circle in $\Si$. Even though
the foliation constructed in \cite{Wa} induces a foliation of $\BH$ by CMC surfaces, since the limit sets in $\Si$ are
nonrectifiable, Wang's results and the results in this paper should be considered in different contexts. Also in
\cite{Wa}, Wang constructs a limit set $\Gamma$ of a quasi-Fuchsian 3-manifold which is similar to the one in
\cite{Ha}, where there cannot be a foliation of $\BH$ by CMC surfaces with asymptotic boundary $\Gamma$. This also
implies the existence of a quasi-Fuchsian $3$-manifold which has no foliation by CMC surfaces.

\end{document}